\documentclass[11pt]{amsart}
\usepackage{amsfonts,amssymb,graphicx}

\input prepictex
\input pictex
\input postpictex

\numberwithin{equation}{section}

\newtheorem{theorem}{Theorem}[section]
\newtheorem{lemma}[theorem]{Lemma}
\newtheorem{proposition}[theorem]{Proposition}
\newtheorem{corollary}[theorem]{Corollary}

\theoremstyle{definition}

\theoremstyle{remark}
\newtheorem{remark}[theorem]{Remark}

\newcommand{\CC}{\mathbb{C}}
\newcommand{\HH}{\mathbb{H}}

\newcommand{\ZZ}{\mathbb{Z}}

\newcommand{\OO}{\mathcal{O}}

\newcommand{\hh}{\mathfrak{h}}

\newcommand{\la}{\langle}
\newcommand{\ra}{\rangle}
\newcommand{\one}{\mathbf{1}}

\newcommand{\ttt}{\mathfrak{t}}
\begin{document}
\title[Combinatorial representation theory of rational Cherednik algebras]{Towards a combinatorial representation theory 
for the rational Cherednik algebra of type $G(r,p,n)$.}
\author{Stephen Griffeth}
\address{Department of Mathematics \\
University of Minnesota \\
Minneapolis, MN 55455 \\}
\email{griffeth@math.umn.edu}

\begin{abstract}
The goal of this paper is to lay the foundations for a combinatorial study, via orthogonal functions and intertwining operators, of category $\OO$ for the rational Cherednik algebra of type $G(r,p,n)$.  As a first application, we give a self-contained and elementary proof of the analog for the groups $G(r,p,n)$, with $r>1$, of Gordon's theorem (previously Haiman's conjecture) on the diagonal coinvariant ring.  We impose no restriction on $p$; the result for $p \neq r$ has been proved by Vale using a technique analogous to Gordon's.  Because of the combinatorial application to Haiman's conjecture, the paper is logically self-contained except for standard facts about complex reflection groups.  The main results should be accessible to mathematicians working in algebraic combinatorics who are unfamiliar with the impressive range of ideas used in Gordon's proof of his theorem.
\end{abstract}

\maketitle

\section{Introduction.}

The purpose of this paper is twofold.  First, we introduce intertwining operators for the rational Cherednik algebra $\HH$ of type $G(r,p,n)$ and carry out enough calculation of the relations they satisfy to be useful for a combinatorial study of the representation theory of $\HH$.  This work forms the basis for the sequels \cite{Gri1}, \cite{Gri2}, and \cite{Gri3}, where we study the combinatorics of the ordinary coinvariant ring and begin the combinatorial study of the lattice of submodules of each standard module for the rational Cherednik algebra.  Our goal is the construction of canonical bases for the standard modules $M(V)$ and their composition factors, including the irreducible quotients $L(V)$.  

Second, we use the intertwining operators to give a new and self-contained (modulo standard facts about complex reflection groups) proof of the analog of Gordon's theorem (see \cite{Gor2} and \cite{Val}) on the diagonal coinvariant ring for $G(r,p,n)$.  The proof here works only when $r>1$ but does not place any restriction on $p$; Vale \cite{Val} establishes the result for $p \neq r$.  We use neither the KZ functor and cyclotomic Hecke algebras (as in \cite{Gor2} and \cite{Val}) nor degeneration from the double affine Hecke algebra (as in \cite{Che}) as there is no DAHA available for the groups $G(r,p,n)$ with $r>2$, although we found the ideas of those papers inspirational.  The existence of this paper also owes much to the foundational papers \cite{DuOp} and \cite{EtGi}.

One of our goals being self-containment, we begin in Section~\ref{PBW Theorem section} with definitions and a sketch of the proof of the Poincar\'e-Birkhoff-Witt theorem for rational Cherednik algebras.  In Section~\ref{diagonal coinvariants} we explain the formalism that connects the rational Cherednik algebra to quotients of the diagonal coinvariant ring of the type conjectured by Haiman \cite{Hai}.  In Section~\ref{Grpn section} we specialize to the case of the groups $G(r,p,n)$ and review the construction given by Dunkl and Opdam in \cite{DuOp} of an important commutative subalgebra of the rational Cherednik algebra $\HH$.  We study the intertwining operators and their basic properties in Section~\ref{intertwiners section} and their action on a particular basis of the polynomial representation in Section~\ref{spectrum section}.  We determine the submodule structure of the polynomial representation of $G(r,p,n)$ in the cases that we will use for the study of the diagonal coinvariant ring in Section~\ref{L section}.  The results of Section~\ref{L section} are similar to, but more detailed than, those contained in \cite{ChEt}.  The extra detail is crucial for the material of Section~\ref{diag for Grpn}, where we study the $\HH$-modules fulfilling the requirements of Section~\ref{diagonal coinvariants} and relevant to the diagonal coinvariant ring.

{\bf Acknowledgements.}  The bulk of this paper is based on a thesis (\cite{Gri}) written at the University of Wisconsin under the direction of Arun Ram.  I am greatly indebted to him for teaching me about rational Cherednik algebras and for suggesting the problems that motivated this work.  I am also grateful to  Drew Armstrong, Ezra Miller, Vic Reiner, and Peter Webb for many interesting discussions during the time this paper was being written. Finally, I would like to thank Richard Vale for correspondence that stimulated an improvement of the results and Iain Gordon for pointing out that Rouquier's results (\cite{Rou}) are relevant to the Bessis-Reiner conjecture on $W$-Catalan numbers.

\section{Definitions and the Poincar\'e-Birkhoff-Witt theorem} \label{PBW Theorem section}

Let $V$ be a finite dimensional vector space over a field $k$, and let $W \subseteq \text{GL}(V)$ be a finite subgroup.  Let $TV$ be the tensor algebra of $V$ and let $kW$ be the group algebra of $W$ over $k$, with basis $t_w$ for $w \in W$ and multiplication $t_w t_v=t_{wv}$.  The \emph{semi-direct product} $TV \rtimes W$ is $TV \otimes_k kW$ with multiplication
\begin{equation}
(f \otimes t_w)( g \otimes t_v)=f(wg) \otimes t_{wv} \quad \hbox{for $f,g \in TV$ and $w,v \in W$.}
\end{equation}  From now on we will drop the tensor signs when it will not cause confusion.  Fix a collection of skew-symmetric forms indexed by the elements of $W$: 
\begin{equation}
\langle \cdot,\cdot,\rangle_w: V \times V \rightarrow k \quad \hbox{for $w \in W$.}
\end{equation}    The \emph{Drinfeld Hecke algebra} $\HH$ corresponding to this data is the quotient of the algebra $TV \otimes_k kW$ by the relations
\begin{equation}
xy-yx=\sum_{w \in W}  \langle x,y \rangle_w t_w \quad \hbox{for $x,y \in V$.}
\end{equation}  

Now let $\hh$ be a finite dimensional $k$-vector space.  A \emph{reflection} is an element $s \in \text{GL}(\hh)$ such that $\text{codim}(\text{fix}(s))=1$.  A \emph{reflection group} is a finite subgroup $W \subseteq \text{GL}(\hh)$ that is generated by the set of reflections it contains.  Assume now that $W$ is a reflection group, let $T$ be the set of reflections it contains, and put $V=\hh^* \oplus \hh$.  Write $\la x,y \ra=x(y)$ for $x \in \hh^*$ and $y \in \hh$.  For each $s \in T$ fix $\alpha_s \in \hh^*$ and $\alpha_s^\vee \in \hh$ with
\begin{equation} \label{alpha eq}
sx=x-\la x,\alpha_s^\vee \ra \alpha_s \quad \hbox{for all $x \in \hh^*$,}
\end{equation} and define skew symmetric bilinear forms $\la \cdot,\cdot \ra_s$ on $V$ by the requirements
\begin{equation} \label{formdef}
\la x, y \ra_s=0 \ \hbox{if $x,y \in \hh$ or $x,y \in \hh^*$ and} \
\la x,y \ra_s=c_s \la \alpha_s,y \ra \la x,\alpha_s^\vee \ra \ \hbox{if $x \in \hh^*$ and $y \in \hh$,}
\end{equation} where $c_s \in k$ satisfy $c_{wsw^{-1}}=c_s$ for all $s \in T$ and $w \in W$.  It is straightforward to check that $\la \cdot,\cdot \ra_s$ does not depend on the choice of $\alpha_s$ and $\alpha_s^\vee$ satisfying \eqref{alpha eq}.  We extend the pairing $\la \cdot,\cdot \ra$ between $\hh$ and $\hh^*$ to a symplectic form on $V$ by requiring that $\hh$ and $\hh^*$ be isotropic.

The \emph{rational Cherednik algebra} corresponding to a reflection group $W$ is the quotient of the semi-direct product $TV \rtimes W$ by the relations
\begin{equation} \label{fundamental relation}
yx-xy=\kappa \la x,y \ra-\sum_{s \in T} \la x,y \ra_s t_s \quad \hbox{for $x,y \in V$,}
\end{equation}  where $\kappa \in k$.  It is thus a special case of the Drinfeld Hecke algebra.  Note that by the definitions of $\la \cdot,\cdot \ra$ and $\la \cdot,\cdot \ra_s$, we have
\begin{equation}
yx-xy=0 \quad \hbox{if $x,y \in \hh$ or $x,y \in \hh^*$,}
\end{equation} and hence there are canonical maps $S(\hh) \rightarrow \HH$ and $S(\hh^*) \rightarrow \HH$.

Now assume that $\HH$ is the Drinfeld Hecke algebra associated to a collection $\la \cdot,\cdot \ra_w$ of skew symmetric forms as above.  We say that \emph{the PBW theorem holds for} $\HH$ if for any basis $x_1,\dots,x_n$ of $V$, the collection $\{x_{i_1} x_{i_2} \dots x_{i_p} t_w | 1 \leq i_1 \leq i_2 \leq \dots \leq i_p \leq n, w \in W \}$ is a basis for $\HH$.  The following theorem was stated in \cite{Dri} and many proofs have now appeared (see, for example, \cite{BaBe}, \cite{EtGi},  \cite{Kha}), and \cite{RaSh}).  Our proof has the virtue of working in arbitrary characteristic; more importantly, it is conceptually extremely simple.  
\begin{theorem}[The Poincar\'e-Birkhoff-Witt theorem for Drinfeld Hecke algebras] \label{PBW}
The PBW theorem holds for $\HH$ if and only if the following two conditions hold:
\item[(a)] $\la vx,vy \ra_{vwv^{-1}}= \la x,y \ra_w$ for all $x,y \in V$ and $v,w \in W$, and
\item[(b)] $\la x,y \ra_w (wz-z)+ \la y,z \ra_w (wx-x)+\la z,x \ra_w (wy-y)=0$ for all $x,y,z \in V$ and $w \in W$.
\end{theorem}
\begin{proof}
First assume the PBW theorem holds for $\HH$.  Let $x,y \in V$ and $v \in W$.  Equating coefficients on both sides of 
\begin{equation}
\sum_{w \in W} \la vx,vy \ra_w t_w=\left[vx,vy \right]=t_v \left[x,y \right] t_v^{-1}=\sum_{w \in W}  \la x,y \ra_w t_{vwv^{-1}}
\end{equation} implies that (a) holds.  Let $x,y,z \in V$.  By the Jacobi identity,
\begin{align*}
0&=\left[ \left[x,y \right],z \right] +\left[\left[y,z\right], x  \right]+\left[ \left[z, x\right] ,y \right]=\left[\sum_{w \in W} \la x,y \ra_w t_w, z \right]+\left[\sum_{w \in W} \la y,z \ra_w t_w, x \right]+\left[ \sum_{w \in W} \la z,x \ra_w t_w,y \right] \\
&=\sum_{w \in W} \left(\la x,y \ra_w(wz-z)+\la y,z \ra_w (wx-x)+\la z,x \ra_w (wy-y) \right) t_w
\end{align*} Now equating coefficients of $t_w$ on both sides implies that (b) holds.

Conversely, assume that (a) and (b) hold.  The defining relations for $\HH$ evidently imply that given any basis $x_1,x_2,\dots,x_n$ of $V$, the set $\{x_{i_1} x_{i_2} \dots x_{i_p} t_w | 1 \leq i_1 \leq i_2 \leq \dots \leq i_p \leq n, w \in W \}$ spans $\HH$.  We will show that these elements are also linearly independent by mimicking the standard proof of the PBW theorem for universal enveloping algebras of Lie algebras: we construct the module that ought to be the left regular representation of $\HH$.  Let $M$ be the vector space with basis consisting of the words $\{x_{i_1} x_{i_2} \dots x_{i_p} t_w | 1 \leq i_1 \leq i_2 \leq \dots \leq i_p \leq n, w \in W \}$.  Define operators $l_x$ and $l_v$ on $M$ for $x \in V$ and $v \in W$ inductively as follows:
\begin{equation}
l_x.t_w=x t_w, \ l_v.t_w=t_{vw},
\end{equation} and for $p \geq 1$,
\begin{equation}
l_{x_i}.x_{i_1}\dots x_{i_p} t_w=\begin{cases} x_i x_{i_1} \dots x_{i_p} t_w \quad \hbox{if $i \leq i_1$,} \\ l_{x_{i_1}}.l_{x_i}.x_{i_2}\dots x_{i_p} + \sum_{v \in W} \la x_i,x_{i_1} \ra_v l_v.x_{i_2} \dots x_{i_p} t_w \quad \hbox{if $i>i_1$,}\end{cases}
\end{equation} and
\begin{equation}
l_v.x_{i_1} \dots x_{i_p} t_w=l_{vx_{i_1}}.l_v.x_{i_2}\dots x_{i_p} t_w.
\end{equation}  A straightforward but lengthy calculation shows that these operators satisfy the defining relations for $\HH$.  It follows that $M$ is an $\HH$-module, with $x$ acting by $l_x$ and $t_w$ acting by $l_w$.  Suppose that there is a relation in $\HH$ of the form
\begin{equation*}
\sum_{\substack{w \in W \\ 1 \leq i_1 \leq \dots \leq i_p \leq n}} a_{i_1 \dots i_p, w} x_{i_1} \dots x_{i_p} t_w=0,
\end{equation*} with $a_{i_1 \dots i_p, w} \in k$.  Applying both sides of this relation to the element $1=t_1 \in M$ implies that all the coefficients $a_{i_1 \dots i_p, w}$ are zero, and the proof is complete.
\end{proof}
\begin{corollary} \label{PBW theorem for RCAs}
Let $\HH$ be the rational Cherednik algebra corresponding to the reflection group $W$.  Then the multiplication map $S(\hh^*) \otimes_k S(\hh) \otimes_k kW \rightarrow \HH$ is an isomorphism.
\end{corollary} 
\begin{proof}
The result follows from the previous theorem once we check that conditions (a) and (b) of that theorem hold.  Condition (a) is straightforward to verify.  Condition (b) holds trivially for $w \notin T$.  If $s \in T$ and $z \in \text{fix}_V(s)$ then the definitions \eqref{alpha eq} and \eqref{formdef} imply $\la z,x \ra_s=0$ for all $x \in V$.  Therefore the radical of $\la \cdot, \cdot \ra_s$ has codimension at most $2$.  If $\la \cdot,\cdot \ra_s \neq 0$ then we may choose $x,y \in V$ with $\la x,y \ra_s=1$,  and for any $z \in V$ we have
\begin{equation}
z=\la x,z \ra_s y-\la y,z \ra_s x +f \quad \text{with} \ f \in \text{fix}_V(s).
\end{equation}  Applying $s-1$ to both sides and rearranging terms shows that the identity (b) holds for $x,y$ and arbitrary $z$.  In general, identity (b) holds trivially if every two of $x,y,z$ are linearly dependent modulo the radical of $\la \cdot,\cdot \ra_s$, so we are reduced to the case just treated.  
\end{proof} 

There is an important filtration of $\HH$ defined by
\begin{equation} \label{filt def}
\HH^{\leq m}=F \text{-span} \{x_{i_1} \cdots x_{i_l} t_w \ | \ l \leq m, \ w \in W, \ \text{and} \ x_{i_j} \in \hh^* \oplus \hh \}.
\end{equation}  By the PBW theorem for $\HH$, the associated graded algebra of $\HH$ with respect to this filtration is the semidirect product $S(\hh^* \oplus \hh) \rtimes W$.

Our next proposition is a fundamental computation.  It expresses some commutators in $\HH$ as linear combinations of derivatives and divided differences of elements of $S(\hh^*)$ and $S(\hh)$.  For $y \in \hh$, we write $\partial_y$ for the derivation of $S(\hh^*)$ determined by
\begin{equation}
\partial_y(x)=\la x,y \ra \quad \hbox{for $x \in \hh^*$,}
\end{equation} and we define a derivation $\partial_x$ of $S(\hh)$ analogously.
\begin{proposition} \label{fund comm rel}
Let $y \in \hh$ and $f \in S(\hh^*)$.  Then
\begin{equation} \label{cf fy}
yf-fy=\kappa \partial_y f - \sum_{s \in T} c_s \la \alpha_s,y\ra \frac{f-sf}{\alpha_s} t_s.
\end{equation}  Similarly, for $x \in \hh^*$ and $g \in S(\hh)$, we have
\begin{equation} \label{cf fx}
gx-xg=\kappa \partial_x g-\sum_{s \in T} c_s \la x,\alpha_s^\vee\ra t_s \frac{g-s^{-1}g}{\alpha_s^\vee}.
\end{equation}
\end{proposition}
\begin{remark}
Note the placement of $t_s$ in the second formula.  In practice, it is sometimes convenient to rewrite it as
\begin{equation}
gx-xg=\kappa \partial_x g-\sum_{s \in T} c_s \la x,\alpha_s^\vee \ra \frac{sg-g}{s\alpha_s^\vee} t_s.
\end{equation}
\end{remark}
\begin{proof}
Observe if $f=x \in \hh^*$, the first formula to be proved is
\begin{equation*}
yx-xy=\kappa \la x,y \ra -\sum_{s \in T} c_s \la \alpha_s,y\ra \frac{x-sx}{\alpha_s} t_s,
\end{equation*} and the right hand side may be rewritten as
\begin{equation*}
\kappa \la x,y \ra-\sum_{s \in T} c_s \la \alpha_s,y\ra \la x,\alpha_s^\vee \ra t_s,
\end{equation*} so that the formula to be proved is one of the defining relations for $\HH$.  We proceed by induction on the degree of $f$.  Assume we have proved the result for $h \in S^d(\hh^*)$ and all $d \leq m$.  For $f,g \in S^{\leq m}(\hh^*)$, and $y \in \hh$, we have
\begin{align*}
[y,fg]&=[y,f]g+f[y,g] \\
&=\left(\kappa \partial_y(f)-\sum_{s \in T} c_s \la \alpha_s,y \ra \frac{f-sf}{\alpha_s} t_s \right) g
+f\left(\kappa \partial_y(g)-\sum_{s \in T} c_s \la \alpha_s,y \ra \frac{g-sg}{\alpha_s} t_s \right)\\
&=\kappa \left(\partial_y(f) g + f \partial_y(g) \right)
-\sum_{s \in T} c_s \la \alpha_s,y\ra \left( \frac{f-sf}{\alpha_s} sg+ f \frac{g-sg}{\alpha_s} \right) t_s \\
&=\kappa \partial_y(fg)-\sum_{s \in T} c_s \la \alpha_s,y \ra \frac{fg-s(fg)}{\alpha_s} t_s.
\end{align*} by using the inductive hypothesis in the second equality, and the Leibniz rule for $\partial_y$ and a \emph{skew} Leibniz rule for the divided differences in the fourth equality.  This proves the first commutator formula, and the proof of the second one is exactly analogous.  
\end{proof}

Let $V$ be a $kW$-module and define a $S(\hh) \otimes_k kW$ action on $V$ by 
\begin{equation}
f.v=f(0)v \quad \text{and} \quad t_w.v=wv \quad \text{for} \quad w \in W, f \in S(\hh).
\end{equation}  The \emph{standard module} corresponding to $V$ is
\begin{equation} \label{standard def}
M(V)=\text{Ind}_{S(\hh) \otimes_k kW}^{\HH} V.
\end{equation}    The PBW theorem shows that $\HH$ is a free $S(\hh) \otimes_k kW$-module, so that the additive functor $V \mapsto M(V)$ is exact and as a $k$-vector space
\begin{equation}
M(V) \simeq S(\hh^*) \otimes_k V.
\end{equation}  In particular when $V=\one$ is the trivial $kW$-module we obtain from Proposition \ref{fund comm rel}
\begin{equation} \label{DunklopsFormula}
M(\one) \simeq S(\hh^*) \quad \text{with} \quad y.f=\kappa \partial_y f-\sum_{s \in  T} c_s \la \alpha_s,y \ra \frac{f-sf}{\alpha_s}
\end{equation} for $y \in \hh$ and $f \in S(\hh^*)$.  These are the famous \emph{Dunkl operators}.  From our point of view, the fact that they commute is a consequence of the PBW theorem, though it is possible to prove the commutativity independently (\cite{DuOp}, for instance) and then use it to establish the PBW theorem.

The definition \eqref{standard def} implies that the module $M(V)$ has the following universal property: given an $\HH$-module $M$ and a $W$-stable subspace $U \subseteq M$ such that $V \cong U$ as $W$ modules and $y.U=0$ for all $y \in \hh$ there is a unique $\HH$-module homomorphism $M(V) \rightarrow M$ which restricts to the given isomorphism $V \cong U$.

Define the element $h \in \HH$ by
\begin{equation}
h=\sum_{i=1}^n x_i y_i + \sum_{s \in T} c_s (1-t_s),
\end{equation} where $x_i$ is a basis of $\hh^*$ and $y_i$ is the dual basis of $\hh$.  Calculations using the defining relations for $\HH$ show that
\begin{equation}
\left[ h,x \right]=\kappa x, \quad \left[h,y \right]=-\kappa y, \quad \text{and} \quad \left[h,t_w \right] =0
\end{equation} for $x \in \hh^*$, $y \in \hh$, and $w \in W$.  Thus if $\kappa=1$ and $V$ is an irreducible $W$-module the $h$ action on the Verma module $M(V)$ is given by
\begin{equation} \label{h grading}
h.fv=(\text{deg}(f)+c_V) fv \quad \hbox{for $f \in S(\hh^*)$ homogeneous and $v \in V$,}
\end{equation} where $c_V$ is the scalar by which $\sum_{s \in T} c_s (1-t_s)$ acts on $V$.

When $\kappa=1$ the formula \eqref{h grading} implies each standard module $M(V)$ has a unique irreducible quotient $L(V)$.  This paper is primarily concerned with the module $L(\one)$ in those cases related to diagonal coinvariants, but the techniques developed will be applied in the sequel \cite{Gri2} to obtain detailed information on the submodule structure of $M(V)$ for more general representations $V$.

\section{Diagonal coinvariants} \label{diagonal coinvariants}

We now describe a situation in which we can relate $L(\one)$ to the diagonal coinvariant ring
\begin{equation}
R=S(\hh^* \oplus \hh)/I \quad \text{where} \quad I=S(\hh^* \oplus \hh)^W_+ S(\hh^* \oplus \hh).
\end{equation}  Recall the filtration \eqref{filt def} of $\HH$ with $\hh^* \oplus \hh$ in degree $1$ and $\CC W$ in degree $0$.  By \eqref{fundamental relation} and the PBW theorem, the associated graded algebra of $\HH$ with respect to this filtration is $S(\hh^* \oplus \hh) \rtimes W$.  It is this fact that was used in \cite{Gor2} (and that we will use in Theorem \ref{diag coinv}) to make the connection to diagonal coinvariants.

In this section, we assume that we are working with an irreducible complex reflection group $W$ of rank $n$.  When $A$ is a graded vector space, we write $A_i$ for the $i$th graded piece. 

\begin{lemma} \label{Koszul lemma}
Assume that there is an irreducible module $V$ of dimension $n$ such that there is an exact sequence $M(V) \rightarrow M(\one) \rightarrow L(\one) \rightarrow 0$ with $L(\one)$ finite dimensional.  Then the Koszul complex
\begin{equation*}
0 \rightarrow S(\hh^*) \otimes \Lambda^n V \rightarrow \cdots \rightarrow S(\hh^*) \otimes \Lambda^1 V \rightarrow S(\hh^*) \rightarrow L(\one) \rightarrow 0
\end{equation*} is exact and the maps are maps of $\HH$-modules.  The graded $W$-character of $L(\one)$ is
\begin{equation*}
\sum_{i \geq 0} \text{tr}(w,L(\one)_i) t^i=\frac{\text{det}(1-t^k w_V)}{\text{det}(1-t w_{\hh^*})}
\end{equation*} where the image of $V$ in $M(\one)$ lies in degree $k$ and $w_V$ and $w_{\hh^*}$ denote $w$ regarded as an endomorphism of $V$ and $\hh^*$, respectively.
\end{lemma}
\begin{proof}
Since $V$ is $n$-dimensional and $L(\one)$ is finite dimensional, the image of $V$ under the map $M(V) \rightarrow M(\one)$ is spanned by a regular sequence.  Hence the Koszul complex is exact.  As a vector space $M(\Lambda^i V) \cong S(\hh^*) \otimes \Lambda^i V$, and using this identification the vector spaces in the Koszul complex are $\HH$-modules.  By assumption the first map is a map of $\HH$-modules; its kernel is therefore an $\HH$-submodule.  The kernel is generated as an $S(\hh^*)$-submodule by $\Lambda^2 V$ and it follows that the second map in the Koszul complex is a map of $\HH$-modules.  One proves in the same way, by induction on $i$, that the $i$th map in the Koszul complex is a map of $\HH$-modules.  That the graded $W$-character of $L(\one)$ is as asserted is a routine calculation using the Koszul resolution.
\end{proof}

If $V$ is an irreducible $W$-module of dimension $l$, its \emph{exponents} are the integers $e_1 \leq e_2 \leq \cdots \leq e_l$ defined by the equation
\begin{equation*}
\sum_{i \geq 0} [(S(\hh^*)/J)_i:V] t^i=\sum_{i=1}^r t^{e_i},
\end{equation*} where $J$ is the ideal generated by the positive degree elements of the invariant ring $S(\hh^*)^W$ and $(S(\hh^*)/J)_i$ denotes the $i$th graded piece.  Since $W$ is a complex reflection group, the invariant ring $S(\hh^*)^W$ is generated by $n$ algebraically independent polynomials $f_1,\dots,f_n$ with degrees $d_1 \leq d_2 \leq \cdots \leq d_n$.

An irreducible representation $V$ is \emph{free} if $(S(\hh^*) \otimes \Lambda^* V^*)^W$ is a free exterior algebra over $S(\hh^*)^W$.  By \cite{OrSo} Theorem 3.1, $V$ is free if $e_1+\cdots+e_l$ is the (unique) exponent of $\Lambda^l V$.  By the proof of \cite{OrSo} Theorem 3.3, the Galois conjugates of the reflection representation $\hh$ are all free.  In Lemma \ref{free reps} we will exhibit some other examples of free representations of the group $G(r,p,n)$ that will be relevant to the diagonal coinvariant ring.  

\begin{theorem} \label{diag coinv}
With assumptions as in Lemma \ref{Koszul lemma}, assume moreover that $V$ is free and the image of $V$ in $M(\one)$ lies in degree $k$ for an integer $k$ such that the multisets $\{k-e_i \}_{i=1}^n$ and $\{d_i \}_{i=1}^n$ are equal.  Then there is a unique occurrence of $\Lambda^n V$ in $L(\one)$, lying in degree $e_1+\cdots +e_n$.  Let $v \in L(\one)$ span this occurrence of $\Lambda^n V$, and filter $L(\one)$ by
\begin{equation*}
L(\one)^{\leq i}=\HH^{\leq i}.v.
\end{equation*}  Then the map $\text{gr} \HH \rightarrow \text{gr} L(\one)$ restricts to a surjection $R \rightarrow \text{gr} L(\one)$, which has $W$-character given by Lemma \ref{Koszul lemma}.  Finally, the image in $\text{gr} L(\one)$ of $S(\hh^*)$ is isomorphic to the ordinary coinvariant ring $S(\hh^*)/J$.  
\end{theorem}
\begin{proof}
In light of Lemma \ref{Koszul lemma}, the occurrences of $\Lambda^n V$ in $L(\one)$ are given by the formula
\begin{equation*}
\sum \left[L(1)_j : \Lambda^n V \right] t^j=\sum (-1)^i \left[ S(\hh^*)_{j-ik} \otimes \Lambda^i V : \Lambda^n V \right] t^j,
\end{equation*}  and we compute the occurrences of $\Lambda^n V$ in $S(\hh^*) \otimes \Lambda^* V$ by use of the $W$-equivariant isomorphism $\Lambda^i V \otimes \Lambda^n V^* \cong \Lambda^{n-i} V^*$.  Thus
\begin{equation*}
\left[S(\hh^*)_j \otimes \Lambda^i V:\Lambda^n V \right]
=\text{dim}_\CC (S(\hh^*)_j \otimes \Lambda^{n-i} V^*)^W.
\end{equation*}  On the other hand, the assumption that $(S(\hh^*) \otimes \Lambda^* V^*)^W$ is a free exterior algebra over $S(\hh^*)^W$ implies that 
\begin{equation*}
\sum \text{dim}_\CC(S(\hh^*)_j \otimes \Lambda^i V^*)^W q^i t^j=\prod_{i=1}^n \frac{1+q t^{e_i}}{1-t^{d_i}}.
\end{equation*}  Thus
\begin{align*}
\sum \left[L(1)_j : \Lambda^n V \right] t^j&=\sum (-1)^i \left[ S(\hh^*)_{j} \otimes \Lambda^i V : \Lambda^n V \right] t^{j+ik} 
=\sum (-1)^i \text{dim}_\CC (S(\hh^*)_j \otimes \Lambda^{n-i} V^*)^W t^{j+ik} \\
&=\sum (-1)^{n-i} \text{dim}_\CC (S(\hh^*)_j \otimes \Lambda^{i} V^*)^W t^{j+(n-i)k} \\
&=(-1)^n t^{nk} \left[ \sum \text{dim}_\CC(S(\hh^*)_j \otimes \Lambda^i V^*)^W q^i t^j \right]_{q=-t^{-k}} \\
&=(-1)^n t^{nk} \prod_{i=1}^n \frac{1- t^{e_i-k}}{1-t^{d_i}}
=t^{e_1+\cdots+e_n} \prod_{i=1}^n \frac{1-t^{k-e_i}}{1-t^{d_i}}=t^{e_1+\cdots+e_n}.
\end{align*}  Upon filtering $L(\one)$ and taking the corresponding associated graded module as in the statement of the theorem, it follows that there is a unique copy of $\Lambda^n V$ in $\text{gr} \HH$.  Since $L(\one)$ is an irreducible $\HH$-module, the map
\begin{equation}
\begin{matrix}
\text{gr} \HH & \longrightarrow & \text{gr} L(\one) \\
f & \longrightarrow & f.v
\end{matrix}
\end{equation} is surjective.  By the PBW-theorem, $\text{gr} \HH=S(\hh^* \oplus \hh) \otimes \CC W$, and the above map remains surjective upon restriction to $S(\hh^* \oplus \hh)$.  Since $v$ is the unique occurrence of $\Lambda^n V$ in $\text{gr} L(\one)$, we have $S(\hh^* \oplus \hh)^W_+.v=0$ and it follows that $\text{gr}L(\one)$ is a quotient of the diagonal coinvariant ring.  

Finally, by Lemma \ref{Koszul lemma} the top degree piece of $L(\one)$ lies in degree 
\begin{equation*}
nk=e_1+\cdots + e_n+d_1+\cdots +d_n=e_1+\cdots+e_n+N
\end{equation*}
and the the socle of the ordinary coinvariant ring $S(\hh^*)/J$ lies in degree $N$.  Since $L(\one)$ is irreducible the map
\begin{equation}
\begin{matrix}
S(\hh^*) \otimes S(\hh) & \longrightarrow & L(\one) \\
f & \longrightarrow & f.v
\end{matrix}
\end{equation} is surjective.  It follows that the socle of $S(\hh^*)/J$ is not in the kernel of the induced map, and hence the coinvariant ring is the image of $S(\hh^*)$.
\end{proof}

If in addition to the hypotheses of Lemma \ref{Koszul lemma} we assume that $V^*$ is free, then an analogous calculation shows that 
\begin{equation}
\sum \text{dim}_\CC L(\one)_i^W t^i =\prod_{i=1}^n \frac{1-t^{k+e_i'}}{1-t^{d_i}}
\end{equation} where $e_1',\dots,e_n'$ are the exponents of $V^*$.  This fact establishes a connection to the conjectural $t$-analog of the $W$-Catalan number discovered by Bessis and Reiner (\cite{BeRe}): if, with the assumptions of Lemma \ref{Koszul lemma}, $W$ is a complex reflection group that can be generated by $n$ reflections, $k=h+1$ where $h=d_n$ is the largest degree (i.e., \emph{Coxeter number}), and $V=\hh^*$ then
\begin{equation}
\sum \text{dim}_\CC L(\one)_i^W t^i =\prod_{i=1}^n \frac{1-t^{h+d_i}}{1-t^{d_i}}.
\end{equation}  It does not seem unreasonable to expect that for most of the exceptional complex reflection groups $W$ and appropriate values of the parameters in the definition of $\HH$, the module $L(\one)$ gives rise to both a nice quotient of the diagonal coinvariant ring and a $t$-analog of the $W$-Catalan number.

\section{The rational Cherednik algebra for $G(r,p,n)$.} \label{Grpn section}

Let $G(r,1,n)$ be the group of $n \times n$ monomial matrices whose entries are $r$th roots of $1$.  Let 
\begin{equation}
\zeta=e^{2 \pi i/r} \quad \text{and} \quad \zeta_i^l=\text{diag}(1,\dots,\zeta^l,\dots,1), \quad \text{for} \quad 1 \leq i \leq n.
\end{equation} Let 
\begin{equation}
s_i=s_{i,i+1}, \quad \text{where} \quad s_{ij}=(ij), \quad \text{for} \quad 1 \leq i < j \leq n,
\end{equation} is the transposition interchanging $i$ and $j$.  There are $r$ conjugacy classes of reflections in $G(r,1,n)$:
(a)  The reflections of order two:
\begin{equation}
\zeta_i^l s_{ij} \zeta_i^{-l}, \quad \text{for} \quad 1 \leq i < j \leq n, \quad 0 \leq l \leq r-1, 
\end{equation} and
(b)  the remaining $r-1$ classes, consisting of diagonal matrices
\begin{equation}
\zeta_i^{l}, \quad \text{for} \quad 1 \leq i \leq n, \quad 1 \leq l \leq r-1,
\end{equation} where $\zeta_i^l$ and $\zeta_j^k$ are conjugate if and only if $k=l$.

Let
\begin{equation*}
y_i=(0,\dots,1,\dots,0)^t \quad \text{and} \quad x_i=(0,\dots,1,\dots,0)
\end{equation*} have $1$'s in the $i$th position and $0$'s elsewhere, so that $y_1,\dots,y_n$ is the standard basis of $\hh=\CC^n$ and $x_1,\dots,x_n$ is the dual basis in $\hh^*$.  If
\begin{equation}
\alpha_s=\zeta^{-l-1}x_i, \quad \alpha_s^\vee=(\zeta^{l+1}-\zeta) y_i, \quad \text{for} \quad s=\zeta_i^l, 
\end{equation} and
\begin{equation}
\alpha_s=x_i-\zeta^l x_j, \quad \alpha_s^\vee=y_i-\zeta^{-l} y_j, \quad \text{for} \quad s=\zeta_i^l s_{ij} \zeta_i^{-l}.
\end{equation} then
\begin{equation*}
sx=x-<x,\alpha_s^\vee> \alpha_s \quad \text{and} \quad s^{-1}(y)=y-<\alpha_s,y>\alpha_s^\vee,
\end{equation*} for $s \in T$, $x \in \hh^*$, and $y \in \hh$.  We relabel the parameters defining $\HH$ by letting
\begin{equation}
c_0=c_{s_1} \quad \text{and} \quad c_i=c_{\zeta_1^i} \quad \text{for} \quad 1 \leq i \leq r-1.
\end{equation}

\begin{proposition}
The rational Cherednik algebra for $W=G(r,1,n)$ with parameters $\kappa,c_0,c_1,\dots,c_{r-1}$ is the algebra generated by $\CC[x_1,\dots,x_n]$, $\CC[y_1,\dots,y_n]$, and $t_w$ for $w \in W$ with relations 
\begin{equation*}
t_w t_v=t_{wv}, \quad t_w x=(wx) t_w, \quad \text{and} \quad t_w y=(wy) t_w,
\end{equation*} for $w,v \in W$, $x \in \hh^*$, and $y \in \hh$, 
\begin{equation} \label{cf2}
y_i x_j=x_j y_i+c_0 \sum_{l=0}^{r-1} \zeta^{-l} t_{\zeta_i^l s_{ij} \zeta_i^{-l}},
\end{equation} for $1\leq i \neq j \leq n$, and 
\begin{equation} \label{cf3}
y_i x_i=x_i y_i+\kappa-\sum_{l=1}^{r-1} c_l (1-\zeta^{-l}) t_{\zeta_i^l} -c_0 \sum_{j \neq i} \sum_{l=0}^{r-1} t_{\zeta_i^l s_{ij} \zeta_i^{-l}},
\end{equation} for $1 \leq i \leq n$.
\end{proposition}
\begin{proof}
This is just a matter of rewriting formula \eqref{fundamental relation} using our $G(r,1,n)$-specific notation.  For $1 \leq i<j \leq n$,
\begin{align*}
y_i x_j&=x_j y_i+\kappa <x_j,y_i> \\
& \qquad -c_0 \sum_{1 \leq k <m \leq n} \sum_{l=0}^{r-1} <x_k-\zeta^l x_m,y_i><x_j,y_k-\zeta^{-l}y_m> t_{\zeta_k^l s_{km} \zeta_k^{-l}} \\
& \qquad -\sum_{k=1}^n \sum_{l=1}^{r-1} c_l <\zeta^{-l-1} x_k,y_i><x_j,(\zeta^{l+1}-\zeta)y_k> t_{\zeta_k^l} \\
&=x_j y_i+\kappa \cdot 0-c_0 \sum_{l=0}^{r-1} (-\zeta^{-l}) t_{\zeta_i^l s_{ij} \zeta_i^{-l}}-0=x_j y_i+c_0 \sum_{l=0}^{r-1} \zeta^{-l} t_{\zeta_i^l s_{ij} \zeta_i^{-l}}.
\end{align*}  The calculation for $1 \leq j < i \leq n$ is similar.  For $i=j$,
\begin{align*}
y_i x_i&=x_i y_i+\kappa<x_i,y_i> \\
& \qquad -c_0 \sum_{1 \leq k<m \leq n} \sum_{l=0}^{r-1} <x_k-\zeta^l x_m,y_i><x_i,y_k-\zeta^{-l} y_m> t_{\zeta_k^l s_{km} \zeta_k^{-l}} \\
& \qquad -\sum_{k=1}^n \sum_{l=1}^{r-1} c_l <\zeta^{-l-1} x_k,y_i><x_i,(\zeta^{l+1}-\zeta)y_k> t_{\zeta_k^l} \\
&=x_i y_i+\kappa-c_0 \sum_{1 \leq i < m \leq n} \sum_{l=0}^{r-1} t_{\zeta_i^l s_{im} \zeta_i^{-l}}
-c_0 \sum_{1 \leq k<i \leq n} \sum_{l=0}^{r-1} t_{\zeta_k^l s_{ik} \zeta_k^{-l}}
-\sum_{l=1}^{r-1} c_l (1-\zeta^{-l}) t_{\zeta_i^l}.
\end{align*}
\end{proof}  Most of the equations that occur later on are simpler in terms of a certain reparametrization.  For $j \in \ZZ$ define
\begin{equation} \label{ds definition}
d_j=\sum_{l=1}^{r-1} \zeta^{lj} c_l.
\end{equation}  It follows that $d_0+d_1+\cdots+d_{r-1}=0$ and that for $1 \leq l \leq r-1$
\begin{equation}
c_l=\frac{1}{r} \sum_{j=0}^{r-1} \zeta^{-lj} d_j.
\end{equation}  The defining relation \eqref{cf3} becomes
\begin{equation} \label{cf3'}
y_i x_i=x_i y_i+\kappa-\sum_{j=0}^{r-1}(d_j-d_{j-1}) \epsilon_{ij}-c_0 \sum_{j \neq i} \sum_{l=0}^{r-1} t_{\zeta_i^l s_{ij} \zeta_i^{-l}},
\end{equation} where for $0 \leq j \leq r-1$, the primitive idempotents for the cyclic reflection subgroup of $W$ generated by $\zeta_i$ are 
\begin{equation} \label{idems def}
\epsilon_{ij}=\frac{1}{r} \sum_{l=0}^{r-1} \zeta^{-lj} t_{\zeta_i^l}.
\end{equation}

The complex reflection group $G(r,p,n)$ is the subgroup of $G(r,1,n)$ consisting of those matrices so that the product of the non-zero entries is an $r/p$th root of $1$.  The reflections in $G(r,p,n)$ are

(a) $\zeta_i^{l} s_{ij} \zeta_i^{-l}$ for $1 \leq i < j \leq n$ and $0 \leq l \leq r-1$, and  

(b) $\zeta_i^{lp}$ for $1 \leq i \leq n$ and $0 \leq l \leq r/p-1$.

\noindent When $n \geq 3$, the rational Cherednik algebra for $G(r,p,n)$ is the subalgebra of the rational Cherednik algebra $\HH$ for $G(r,1,n)$ with parameters
$$c_l=0 \ \hbox{if $p$ does not divide $l$},$$ generated by $\CC[x_1,\dots,x_n]$, $\CC[y_1,\dots,y_n]$, and $\CC G(r,p,n)$.  Although this is not strictly speaking true for $n=2$, our results on the diagonal coinvariant ring still go through in that case except when $p=r=2$.

Although not strictly necessary for the results of this paper, it seems worthwhile to mention here that when $c_l=0$ for $l$ not divisible by $p$ there is a cyclic group of automorphisms, generated by
\begin{equation}
x \mapsto x, \quad y \mapsto y, \quad s_i \mapsto s_i, \quad t_{\zeta_j} \mapsto \zeta^{r/p} t_{\zeta_j},
\end{equation} for $x \in \hh^*$, $y \in \hh$, $1 \leq i \leq n-1$, and $1 \leq j \leq n$, of the rational Cherednik algebra $\HH$ for $G(r,1,n)$ so that the rational Cherednik algebra for $G(r,p,n)$ is the fixed subalgebra.  The version of Clifford theory given in \cite{RaRa} therefore applies to deduce representation theoretic results for the $G(r,p,n)$ RCA from those for the $G(r,1,n)$ RCA.

From now on, $p$ dividing $r$ will be fixed and we work with the rational Cherednik algebra $\HH$ for $G(r,p,n)$.  Note that with the parameters $c_l=0$ for $l$ not divisible by $p$, we have $d_j=d_k$ if $j=k$ mod $r/p$.

Our first goal is identify a certain commutative subalgebra $\ttt$ of $\HH$.  Later on we will use the subalgebra $\ttt$  to diagonalize the standard module $M(\one)$ (as in \cite{DuOp}), a result which is generalized in the paper \cite{Gri2}.  For $1 \leq i \leq n$ define
\begin{equation} \label{z def}
z_i=y_i x_i+c_0 \phi_i \quad \text{for} \quad 1 \leq i \leq n, \quad \text{where} \quad \phi_i=\sum_{1 \leq j <i} \sum_{l=0}^{r-1} t_{\zeta_i^l s_{ij} \zeta_i^{-l}}.
\end{equation}  The following proposition is proved in \cite{DuOp}.  
\begin{proposition} \label{zs comm}
The elements $z_1,\dots,z_n$ of $\HH$ are pairwise commutative:
\begin{equation*}
z_i z_j=z_j z_i \quad \text{for} \quad 1 \leq i,j \leq n.
\end{equation*}
\end{proposition}
\begin{proof}
We begin by computing
\begin{align*}
\left[ y_i x_i,y_j x_j \right]&=y_i x_i y_j x_j-y_j x_j y_i x_i 
=y_i(x_i y_j-y_j x_i)x_j+y_j(y_i x_j-x_j y_i) x_i \\
&=-y_i\left( c_0 \sum_{l=0}^{r-1} \zeta^{-l} t_{\zeta_j^l s_{ij} \zeta_j^{-l}} \right) x_j+y_j \left( c_0 \sum_{l=0}^{r-1} \zeta^{-l} t_{\zeta_i^l s_{ij} \zeta_i^{-l}} \right) x_i \\
&=-y_i x_i \left( c_0 \sum_{l=0}^{r-1} t_{\zeta_j^l s_{ij} \zeta_j^{-l}} \right) + \left( c_0 \sum_{l=0}^{r-1} t_{\zeta_i^{l} s_{ij} \zeta_i^{-l}} \right) y_i x_i 
=-\left[y_i x_i, c_0 \sum_{l=0}^{r-1} t_{\zeta_i^l s_{ij} \zeta_i^{-l}} \right].
\end{align*}  Thus 
\begin{equation} \label{cf}
\left[ y_i x_i,y_j x_j+c_0 \sum_{l=0}^{r-1} t_{\zeta_i^l s_{ij} \zeta_i^{-l}} \right]=0
\end{equation}  Let
\begin{equation*}
\psi_i=\phi_1+\dots+\phi_i=\sum_{1 \leq j < k \leq i} \sum_{l=0}^{r-1} t_{\zeta_k^l s_{jk} \zeta_k^{-l}}.
\end{equation*}  Then $\psi_i$ is a conjugacy class sum and therefore a central element of the group algebra of $G(r,1,i)$.  It follows that $\psi_i$ commutes with $\psi_1,\dots,\psi_i$.  Therefore $\psi_1,\psi_2,\dots,\psi_n$ are pairwise commutative and hence $\phi_1,\dots,\phi_n$ are pairwise commutative.

Using the commutativity of the $\phi_i$, the commutator formula \eqref{cf}, and the fact that $y_j x_j$ commutes with $\phi_i$ for $i < j$, we assume $i<j$ and compute
\begin{align*}
[z_i,z_j]&=[y_i x_i+c_0 \phi_i, y_j x_j+c_0 \phi_j]=[y_i x_i, y_j x_j + c_0 \phi_j ] \\
&=\left[ y_i x_i, y_j x_j+c_0 \sum_{l=0}^{r-1} t_{\zeta_i^l s_{ij} \zeta_i^{-l}} + c_0 \sum_{1 \leq k \neq i < j} \sum_{l=0}^{r-1} t_{\zeta_j^l s_{jk} \zeta_j^{-l}} \right]=0.
\end{align*}
\end{proof} 

As observed in \cite{Dez}, Proposition 1.1, the relations in the following lemma imply that the subalgebra of $\HH$ generated by $G(r,p,n)$ and $z_1,\dots,z_n$ is isomorphic to the graded Hecke algebra for $G(r,p,n)$ defined in \cite{RaSh} Section 5 (the elements $z_1,\dots,z_n$ are algebraically independent over $\CC$ by the PBW theorem).  
\begin{proposition} \label{comm rels}  Working in the rational Cherednik algebra $\HH$ for $G(r,1,n)$,
\begin{equation} \label{zt1}
z_i t_{\zeta_j} = t_{\zeta_j} z_i \quad \text{for} \quad 1 \leq i,j \leq n,
\end{equation}
\begin{equation} \label{intertwiner1}
z_i t_{s_i}=t_{s_i} z_{i+1} -c_0 \sum_{l=0}^{r-1} t_{\zeta_i^l \zeta_{i+1}^{-l}} \quad \text{for} \quad 1 \leq i \leq n,
\end{equation} and
\begin{equation} \label{intertwiner2}
z_i t_{s_j}=t_{s_j} z_i \quad \text{for} \quad 1 \leq i \leq n \quad \text{and} \quad j \neq i, i+1.
\end{equation}
\end{proposition}
\begin{proof}
First we observe that the elements $t_{\zeta_i}$ and $\phi_j$ commute for all $1 \leq i,j \leq n$.  This is clear if $i > j$; if $i = j$ then
\begin{equation*}
t_{\zeta_j} \phi_j t_{\zeta_j}^{-1}=t_{\zeta_j} \sum_{1 \leq k < j} \sum_{l=0}^{r-1} t_{\zeta_j^l s_{jk} \zeta_j^{-l}} t_{\zeta_j}^{-1}=\sum_{1 \leq k < j} \sum_{l=0}^{r-1} t_{\zeta_j^{l+1} s_{jk} \zeta_j^{-l-1}}=\phi_j;
\end{equation*} a similar computation handles the case $i < j$.  Then \eqref{zt1} follows from
\begin{equation*}
t_{\zeta_i} y_i x_i= \zeta y_i t_{\zeta_i} x_i=\zeta \zeta^{-1} y_i x_i t_{\zeta_i} = y_i x_i t_{\zeta_i}, \quad \text{for} \quad 1 \leq i,j \leq n.
\end{equation*}  For \eqref{intertwiner1},
\begin{align*}
z_i t_{s_i}&=\left( y_i x_i+c_0 \sum_{1 \leq j < i} \sum_{l=0}^{r-1} t_{\zeta_i^l s_{ij} \zeta_i^{-l}} \right) t_{s_i}=t_{s_i} \left(y_{i+1} x_{i+1}+c_0 \sum_{1 \leq j <i} \sum_{l=0}^{r-1} t_{\zeta_{i+1}^l s_{i+1,j} \zeta_{i+1}^{-l}} \right) \\
&=t_{s_i} z_{i+1}-t_{s_i} c_0 \sum_{l=0}^{r-1} t_{\zeta_{i+1}^l s_{i+1,i} \zeta_{i+1}^{-l}}=t_{s_i} z_{i+1}-c_0 \sum_{l=0}^{r-1} t_{\zeta_i^l \zeta_{i+1}^{-l}}.
\end{align*}  Finally, we observe that if $j \neq i, i+1$ then $t_{s_j}$ commutes with $\phi_i$ and with $y_i x_i$, and hence with $z_i=y_i x_i + c_0 \phi_i$.
\end{proof}  

Let 
\begin{equation}
\ttt=\CC[z_1,\dots,z_n,t_{\zeta_1 \zeta_2^{-1}},\dots,t_{\zeta_{n-1}\zeta_n^{-1}},t_{\zeta_1^p},\dots,t_{\zeta_n^p}].
\end{equation}  By Proposition \ref{zs comm} and Lemma \ref{comm rels}, the subalgebra $\ttt$ is a commutative subalgebra of $\HH$.  Our goal is to use $\ttt$ in much the same way as a Cartan subalgebra of a semisimple Lie algebra.

\section{Intertwiners.} \label{intertwiners section}

In this section we will prove many formulas using elements of $G(r,1,n)$ that are not in $G(r,p,n)$; since the rational Cherednik algebra for $G(r,p,n)$ is a subalgebra of a specialization of that for $G(r,1,n)$, these formulas have consequences in the rational Cherednik algebra for $G(r,p,n)$.

The following lemma is a generalization of \eqref{intertwiner1} and \eqref{intertwiner2}.  Let
\begin{equation}
\pi_i=\sum_{l=0}^{r-1} t_{\zeta_i^l \zeta_{i+1}^{-l}}.
\end{equation}
\begin{lemma} \label{zt cf1}
Let $f$ be a rational function of $z_1,\dots,z_n$.  Then
\begin{equation} \label{intertwiner3}
t_{s_i} f=(s_i f) t_{s_i}-c_0 \pi_i \frac{f-s_i f}{z_i - z_{i+1}}, \quad \text{for} \quad 1 \leq i \leq n-1.
\end{equation}
\end{lemma}
\begin{proof}
Observe that if $f$ is $z_i$, $z_{i+1}$, or $z_j$ for $j \neq i,i+1$, then the relation to be proved is \eqref{intertwiner1} and \eqref{intertwiner2}.  Assume the relation \eqref{intertwiner3} is true for rational functions $f$ and $g$.  Then it is evidently true for $f+g$ and $af$ for all $a \in \CC$, and we compute
\begin{align*}
t_{s_i} fg&=\left(t_{s_i} f-(s_i f) t_{s_i} \right) g+ (s_i f) \left(t_{s_i} g-(s_i g) t_{s_i} \right) + (s_i fg) t_{s_i} \\
&=\left( -c_0 \pi_i \frac{f-s_if}{z_i-z_{i+1}} \right) g + (s_i f) \left( -c_0 \pi_i \frac{g-s_ig}{z_i-z_{i+1}} \right) + (s_i fg) t_{s_i} \\
&=(s_i fg) t_{s_i} -c_0 \pi_i \frac{fg-s_i fg}{z_i-z_{i+1}},
\end{align*} so \eqref{intertwiner3} is true for $fg$.  Assuming it is true for the rational function $f$, we compute
\begin{align*}
(t_{s_i} 1/f- (1/s_if) t_{s_i} ) f (s_i f)&= t_{s_i} s_i f- 1/s_i f \left( (s_i f) t_{s_i} -c_0 \pi_i \frac{f-s_i f}{z_i-z_{i+1}} \right) s_i f \\
&=c_0 \pi_i \frac{f-s_i f}{z_i-z_{i+1}},
\end{align*} and dividing by $f (s_i f)$ proves that the relation holds for $1/f$.  Since it holds for $z_1,\dots,z_n$, it is true for all rational functions in $z_1,\dots,z_n$.
\end{proof}
The \emph{intertwining operators} $\sigma_i$ for $1 \leq i \leq n-1$ are
\begin{equation} \label{sigmas def}
\sigma_i=t_{s_i}+\frac{c_0}{z_i-z_{i+1}} \pi_i, \quad \text{where} \quad \pi_i=\sum_{l=0}^{r-1} t_{\zeta_i^l \zeta_{i+1}^{-l}},
\end{equation} and we define intertwining operators $\Phi$ and $\Psi$ by
\begin{equation} \label{phipsidef}
\Phi=x_n t_{s_{n-1} s_{n-2} \dots s_1} \quad \text{and} \quad \Psi=y_1 t_{s_1 s_2 \dots s_{n-1}}.
\end{equation}   The intertwiner $\Phi$ was first defined in Section 4 of \cite{KnSa} where it is used for the symmetric group case.  The other intertwiners were defined for the first time in the author's thesis \cite{Gri}, where some of the results that follow were also recorded.

The intertwiner $\sigma_i$ is well-defined when $\pi_i=0$ or $z_i-z_{i+1} \neq 0$.  The intertwiners are important because, as Lemma \ref{s int} shows, they permute the $z_i$'s from \eqref{z def}.  Our first task is to compute the squares $\sigma_i^2$ of the intertwiners and the products $\Phi \Psi$ and $\Psi \Phi$. Since these compositions all lie in $\ttt$ this calculation is useful for deciding when the intertwiners applied to a $\ttt$-eigenvector (or generalized eigenvector) are non-zero.

\begin{lemma} \label{int squares}
\item[(a)]  For $1 \leq i \leq n-1$, 
\begin{equation*}
\sigma_i^2=1-\left(\frac{c_0 \pi_i}{z_i-z_{i+1}}\right)^2.
\end{equation*}
\item[(b)]
\begin{equation*}
\Psi \Phi=z_1 \quad \text{and} \quad \Phi \Psi=z_n-\kappa+\sum_{j=0}^{r-1} (d_j-d_{j-1}) \epsilon_{1j}.
\end{equation*}
\end{lemma}
\begin{proof}
Using Lemma \ref{zt cf1},
\begin{align*}
\sigma_i^2&=\left(t_{s_i}+\frac{c_0 \pi_i}{z_i-z_{i+1}} \right)\left(t_{s_i}+\frac{c_0 \pi_i}{z_i-z_{i+1}}\right) \\
&=1+t_{s_i} \frac{c_0 \pi_i}{z_i-z_{i+1}}+\frac{c_0 \pi_i}{z_i-z_{i+1}}t_{s_i} + \left(\frac{c_0 \pi_i}{z_i-z_{i+1}} \right)^2 \\
&=1+\frac{c_0 \pi_i}{z_{i+1}-z_i} t_{s_i}- c_0 \pi_i \frac{ \frac{c_0 \pi_i}{z_i-z_{i+1}}-\frac{c_0 \pi_i}{z_{i+1}-z_i} }{z_i-z_{i+1}} + \frac{c_0 \pi_i}{z_i-z_{i+1}}t_{s_i} + \left(\frac{c_0 \pi_i}{z_i-z_{i+1}} \right)^2 \\
&=1-\left(\frac{c_0 \pi_i}{z_i-z_{i+1}} \right)^2.
\end{align*}  This proves (a), and (b) follows from the definition \eqref{phipsidef} and the relation \eqref{cf3'}.
\end{proof}

We define a symmetric group action on $\ttt$ by
letting $S_n$ simultaneously permute $z_1,z_2,\dots,z_n$ and $t_{\zeta_1},\dots,t_{\zeta_n}$.  We also define an automorphism $\phi$ of $\ttt$ by
\begin{equation}
\phi(t_{\zeta_i})=t_{\zeta_{i+1}} \hbox{for $1 \leq i \leq n-1$,} \ \phi(t_{\zeta_n})=\zeta^{-1} t_{\zeta_1}
\end{equation} and
\begin{equation}
\phi(z_i)=z_{i+1} \quad \hbox{for $1 \leq i \leq n-1$ and} \ \phi(z_n)=z_1+\kappa-\sum_{j=0}^{r-1}(d_{j-1}-d_{j-2}) \epsilon_{1j}
\end{equation} where as in \eqref{idems def} $\epsilon_{1j}$ are the primitive idempotents for the cyclic reflection subgroup generated by $\zeta_1$.  
\begin{lemma} \label{s int}
\item[(a)]  For $1 \leq i \leq n-1$ and $f \in \ttt$,
\begin{equation*}
\sigma_i f=(s_i.f) \sigma_i.
\end{equation*} 
\item[(b)]  For $f \in \ttt$,
\begin{equation*}
f \Phi=\Phi (\phi.f) \quad \text{and} \quad f \Psi=\Psi (\phi^{-1}.f)
\end{equation*}
\end{lemma}
\begin{proof}
The commutation relation \eqref{intertwiner1} for $z_i$ and $t_{s_i}$ gives
\begin{align*}
z_i \sigma_i&=z_i \left( t_{s_i}+\frac{c_0 \pi_i}{z_i-z_{i+1}} \right)=t_{s_i} z_{i+1}-c_0 \sum_{l=0}^{r-1} t_{\zeta_i^l \zeta_{i+1}^{-l}}+\frac{c_0 \pi_i z_i}{z_i-z_{i+1}} \\
&=\sigma_i z_{i+1}-\frac{c_0 \pi_i z_{i+1}}{z_i-z_{i+1}}-c_0 \pi_i+\frac{c_0 \pi_i z_i}{z_i-z_{i+1}}=\sigma_i z_{i+1}.
\end{align*}  The proof that $z_{i+1} \sigma_i=\sigma_i z_i$ is exactly analogous, and the fact that $\sigma_i$ and $z_j$ commute if $j \neq i, i+1$ is obvious.

Using the relation $t_{\zeta_i} \pi_i=\pi_i t_{\zeta_{i+1}}$,
\begin{align*}
t_{\zeta_i} \sigma_i&=t_{\zeta_i} \left( t_{s_i}+\frac{c_0 \pi_i}{z_i-z_{i+1}} \right)=t_{s_i} t_{\zeta_{i+1}}+\frac{c_0 \pi_i}{z_i-z_{i+1}} t_{\zeta_{i+1}}=\sigma_i t_{\zeta_{i+1}}.
\end{align*}  The proof that $t_{\zeta_{i+1}} \sigma_i=\sigma_i t_{\zeta_i}$ is the same, and the fact that $\sigma_i$ and $t_{\zeta_j}$ commute if $j \neq i, i+1$ is obvious.  This proves (a).  

Using the commutation formula \eqref{cf3'} for $y_n$ and $x_n$,
\begin{align*}
y_n x_n \Phi&=\left(x_n y_n+\kappa-\sum_{j=0}^{r-1} (d_j-d_{j-1}) \epsilon_{nj}-c_0 \phi_n \right) x_n t_{s_{n-1}} \dots t_{s_1} \\
&=\Phi y_1 x_1+\kappa \Phi- \Phi \sum_{j=0}^{r-1} (d_j-d_{j-1}) \epsilon_{1,j+1} - c_0 \phi_n \Phi.
\end{align*} Hence
\begin{equation*}
z_n \Phi=(y_n x_n + c_0 \phi_n) \Phi= \Phi \left(z_1+\kappa-\sum_{j=0}^{r-1} (d_{j-1}-d_{j-2}) \epsilon_{1j} \right).
\end{equation*}  Let $1 \leq i < n$.  Since
\begin{align*}
y_i x_i \Phi&=y_i x_i x_n t_{s_{n-1} \dots s_1}=\left(x_n y_i + c_0 \sum_{l=0}^{r-1} \zeta^{-l} t_{\zeta_i^l s_{in} \zeta_i^{-l}} \right) x_i t_{s_{n-1} \dots s_1} \\
&=\Phi y_{i+1} x_{i+1} + \Phi c_0 \sum_{l=0}^{r-1} t_{\zeta_{i+1}^l s_{i+1,1} \zeta_{i+1}^{-l}}
\end{align*} and
\begin{align*}
\phi_i \Phi&=\sum_{1 \leq j < i} \sum_{l=0}^{r-1} t_{\zeta_i^l s_{ij} \zeta_i^{-l}} x_n t_{s_{n-1} \dots s_1}=x_n t_{s_{n-1} \dots s_1} \sum_{1 \leq j < i} \sum_{l=0}^{r-1} t_{\zeta_{i+1}^l s_{i+1,j+1} \zeta_{i+1}^{-l}} \\
&=\Phi \left(\phi_{i+1}-\sum_{l=0}^{r-1} t_{\zeta_{i+1}^l s_{i+1,1} \zeta_{i+1}^{-l}} \right)
\end{align*} it follows that
\begin{align*}
z_i \Phi&=(y_i x_i + c_0 \phi_i) \Phi \\
&=\Phi y_{i+1} x_{i+1} +\Phi c_0 \sum_{l=0}^{r-1} t_{\zeta_{i+1}^l s_{i+1,1} \zeta_{i+1}^{-l}}+\Phi c_0 \left(\phi_{i+1}-\sum_{l=0}^{r-1} t_{\zeta_{i+1}^l s_{i+1,1} \zeta_{i+1}^{-l}} \right)=\Phi z_{i+1}.
\end{align*} Finally, 
\begin{align*}
t_{\zeta_i} \Phi&=t_{\zeta_i} x_n t_{s_{n-1} \dots s_1}=x_n t_{s_{n-1} \dots s_1} t_{\zeta_{i+1}}=\Phi t_{\zeta_{i+1}}, \quad \text{and} \\
t_{\zeta_n} \Phi&=t_{\zeta_n} x_n t_{s_{n-1} \dots s_1}=x_n t_{s_{n-1} \dots s_1} \zeta^{-1} t_{\zeta_1}=\Phi(\zeta^{-1} t_{\zeta_1}),
\end{align*} for $1 \leq i < n$.  This proves the formula involving $\Phi$.  The formula for $\Psi$ follows from that for $\Phi$ by using the relations in part (b) of Lemma \ref{int squares}.
\end{proof}

\section{An eigenbasis of $M(\one)$} \label{spectrum section}

Let $\CC[x_1,x_2,\dots,x_n]=M(\one)$ be the polynomial representation of $\HH$.  We will show that for generic choices of the parameters $\kappa$ and $c_i$, the ring $\CC[x_1,\dots,x_n]$ has an $\ttt$-eigenbasis indexed by the set $\ZZ_{\geq 0}^n$ and we will describe how the intertwining operators act on this basis.  

For $\mu \in \ZZ_{\geq 0}^n$, let $v_\mu$ be the maximal length permutation such that
\begin{equation}
v_\mu.\mu=\mu_-, \quad \hbox{where $\mu_-$ is the non-decreasing (anti-partition) rearrangement of $\mu$.}
\end{equation}  We write $\mu_+$ for the partition rearrangement of $\mu$, and define a partial order on $\ZZ_{\geq 0}^n$ by
\begin{equation} \label{order def}
\lambda < \mu \quad \iff \quad \lambda_+ <_d \mu_+ \quad \text{or} \quad \lambda_+=\mu_+ \quad \text{and} \quad v_\lambda < v_\mu,
\end{equation}  where we use the Bruhat order on $S_n$, and $<_d$ denotes dominance order on $\ZZ_{\geq 0}^n$, given by
\begin{equation}
\lambda \leq_d \mu \quad \text{if} \quad \mu-\lambda \in \sum_{i=1}^{n-1} \ZZ_{\geq 0} (\epsilon_i-\epsilon_{i+1}).
\end{equation}  If $\mu_i>\mu_{i+1}$ then
\begin{equation} \label{order fact}
\mu> s_i.\mu+k(\epsilon_i-\epsilon_{i+1}) \quad \hbox{for $0 \leq k <\mu_i-\mu_{i+1}$.}
\end{equation}

The next theorem is the analogue of \cite{Opd} 2.6 in our setting.  It shows that the $z_i$'s are upper triangular as operators on $\CC[x_1,\dots,x_n]$ with respect to the order on $\ZZ_{\geq 0}^n$ defined in \eqref{order def}.  Equivalent results are proved \cite{DuOp} by reduction to the symmetric group case.  The proof we give below is generalized to the modules $M(V)$ for all irreducible $\CC W$-modules $V$ in \cite{Gri2}.

\begin{theorem} \label{Jack construction}
\begin{enumerate}
\item[(a)]  The actions of $t_{\zeta_i}^p$, $t_{\zeta_i^{-1} \zeta_{i+1}}$,  and $z_i$ on $M(\one)$ are given by
\begin{equation*}
t_{\zeta_i}^p.x^\mu=\zeta^{-p\mu_i} x^\mu, \quad t_{\zeta_i^{-1} \zeta_{i+1}}.x^\mu=\zeta^{\mu_i-\mu_{i+1}} x^\mu,
\end{equation*}
and, with $d_j$ as in \eqref{ds definition},
\begin{equation*} \label{z spectrum}
z_i.x^\mu=\left( \kappa(\mu_i+1)-(d_0-d_{-\mu_i-1})-r (v_\mu(i)-1) c_0 \right) x^\mu +\sum_{\nu < \mu} c_{\nu} x^\nu.
\end{equation*}
\item[(b)]  Assuming that the parameters are generic, so $F=\CC(\kappa,c_0,d_1,\dots,d_{r/p-1})$, for each $\mu \in \ZZ_{\geq 0}^n$ there exists a unique $\ttt$-eigenvector $f_{\mu} \in M(\one)$ such that
\begin{equation*}
f_{\mu,T}=x^\mu+\text{lower terms}.
\end{equation*}  
\end{enumerate}
The $\ttt$-eigenvalue of $f_{\mu}$ is determined by the formulas in part (a).
\end{theorem}
\begin{proof}
The statements about the action of $t_{\zeta_i^p}$ and $t_{\zeta_i \zeta_{i+1}^{-1}}$ follow from the commutation relation in the definition of the rational Cherednik algebra and the definition of the representation $M(\one)$.  Using the commutation formula in Proposition \ref{fund comm rel} for $f \in \CC[x_1,\dots,x_n]$ and $y \in \hh$ and the geometric series formula to evaluate the divided differences, we obtain the following formula for the action of $y_i x_i$ on $x^\mu$:
\begin{align*}
y_i.x^{\mu+\epsilon_i}&=\kappa (\mu_i+1) x^{\mu}-c_0 \sum_{1 \leq j<k \leq n} \sum_{l=0}^{r-1} \la x_j-\zeta^{l}x_k,y_i \ra \frac{x^{\mu+\epsilon_i}-\zeta_j^l s_{jk} \zeta_j^{-l} x^{\mu+\epsilon_i}}{x_j-\zeta^l x_k}   \\
& \qquad -\sum_{1 \leq j \leq n} \sum_{l=1}^{r/p-1} c_{lp} \la x_j,y_i \ra \frac{x^{\mu+\epsilon_i}-\zeta_j^{lp} x^{\mu+\epsilon_i}}{ x_j }  \\
&=\kappa(\mu_i+1) x^{\mu}-c_0 \sum_{j \neq i} \sum_{l=0}^{r-1}  \frac{x^{\mu+\epsilon_i}-\zeta_i^l s_{ij} \zeta_i^{-l} x^{\mu+\epsilon_i}}{x_i-\zeta^l x_j}
 -\sum_{l=1}^{r/p-1} c_{lp} (1-\zeta^{-lp(\mu_i+1)}) \zeta^{lp b_i} x^{\mu} \\
 &=(\kappa(\mu_i+1)-\sum_{l=1}^{r/p-1} c_{lp} (1-\zeta^{-lp(\mu_i+1)}) \zeta^{lpb_i}) x^{\mu} 
 - c_0 \sum_{\substack{j \neq i \\ \mu_i \geq \mu_j }} (x^\mu+ \zeta^l x^{\mu+(\epsilon_j-\epsilon_i)}+ \cdots +\zeta^{l(\mu_i-\mu_j)} x^{ s_{ij} \mu} ) \\ 
 &+c_0 \sum_{\substack{j \neq i \\ \mu_j > \mu_i }} \sum_{k=1}^{\mu_j-\mu_i-1} \zeta^{-lk} x^{\mu+k(\epsilon_i-\epsilon_j)} .
\end{align*}

Using this equation and \eqref{order fact} to identify lower terms,
\begin{align*}
z_i .x^\mu&=\left(y_i x_i +c_0 \sum_{1 \leq j <i} \sum_{0 \leq l \leq r-1} t_{\zeta_i^l s_{ij} \zeta_i^{-l}} \right) x^\mu 
=y_i. x^{\mu+\epsilon_i} +c_0 \sum_{1 \leq j < i} \sum_{l=0}^{r-1} t_{\zeta_i^{l} s_{ij} \zeta_i^{-l}}. x^\mu  \\
&= (\kappa(\mu_i+1)-\sum_{l=1}^{r/p-1} c_{lp} (1-\zeta^{-lp(\mu_i+1)}) \zeta^{lp b_i}) x^{\mu}
 - c_0 \sum_{\substack{1 \leq j < i \\ \mu_j \leq \mu_i \\ 0 \leq l \leq r-1}} \zeta^{l(\mu_i-\mu_j)} x^{ s_{ij} \mu}  \\
 & - c_0 \sum_{\substack{ 1 \leq j < i \\ \mu_j < \mu_i \\ 0 \leq l \leq r-1}} x^\mu
 - c_0 \sum_{\substack{ i < j \leq n \\ \mu_j \leq \mu_i \\ 0 \leq l \leq r-1}}x^\mu
 +c_0  \sum_{\substack{1 \leq j < i \\ 0 \leq l \leq r-1}} \zeta^{l(\mu_i-\mu_j)} x^{ s_{ij} \mu}+ \text{lower terms} \\
 &= \Big( \kappa(\mu_i+1)-\sum_{l=1}^{r/p-1} c_{lp} (1-\zeta^{-lp(\mu_i+1)}) \zeta^{lpb_i} -c_0 r
 (v_\mu(i)-1) \Big) x^{\mu}
 + \text{lower terms},
\end{align*}  where to obtain the last line we used the formula
\begin{equation*}
v_\mu(i)=|\{j<i \ | \ \mu_j < \mu_i \}|+|\{j > i \ | \ \mu_j \leq \mu_i \}| + 1.
\end{equation*} Now rewriting things in terms of the $d_j$'s from \eqref{ds definition} proves part (a) of the theorem.

For part (b), simply observe that the coefficient of $\kappa$ in the formula for the action of $z_i$ on $x^\mu$ is $\mu_{i+1}$; it follows that the $\ttt$-eigenspaces are all one-dimensional and hence a simultaneous eigenbasis exists.
\end{proof}

Define the \emph{weight} $\text{wt}(\mu)$of $\mu \in \ZZ_{\geq 0}^n$ to be the $\ttt$-homomorphism mapping $z_i$ to
\begin{equation*}
\kappa(\mu_i+1)-(d_0-d_{-\mu_i-1})-r(v_\mu(i)-1) c_0
\end{equation*} and $\zeta_i$ to $\zeta^{-\mu_i}$.  

\begin{lemma} \label{action lemma}
The action of the intertwiners on the basis $f_\mu$ is given by:
\begin{equation} \
\sigma_i.f_\mu=f_{s_i.\mu} \quad \hbox{if $\mu_i<\mu_{i+1}$ or $\mu_i \neq \mu_{i+1}$ mod $r$,}
\end{equation}
\begin{equation}
\sigma_i.f_\mu=0 \quad \hbox{if $\mu_i=\mu_{i+1}$,}
\end{equation}
\begin{equation}
\sigma_i.f_\mu=\frac{(\delta-rc_0)(\delta+rc_0)}{\delta^2}  f_{s_i.\mu} \quad \hbox{if $\mu_i=\mu_{i+1}$ mod $r$ and $\mu_i > \mu_{i+1}$,}
\end{equation} where
\begin{equation*}
\delta=\kappa(\mu_{i}-\mu_{i+1})-c_0 r (w_+(\mu)^{-1}(i+1)-w_+(\mu)^{-1}(i)),
\end{equation*}
\begin{equation}
\Phi.f_\mu=f_{\phi.\mu},
\end{equation}
\begin{equation}
\Psi.f_\mu=0 \quad \hbox{if $\mu_n=0$,}
\end{equation} and
\begin{equation}
\Psi.f_\mu=\left(\kappa \mu_n-(d_0-d_{-\mu_n})-c_0 r (v_\mu(n)-1)  \right) f_{\psi.\mu} \quad \hbox{if $\mu_n \neq 0$.}
\end{equation}
\end{lemma}
\begin{proof}
We will establish the formulas for $\sigma_i$; the formulas for $\Phi$ and $\Psi$ are proved in an analogous fashion.  If $\mu_i < \mu_{i+1}$, then for all $\nu \leq \mu$ one has $s_i.\nu \leq s_i.\mu$ and it follows that the leading term of $\sigma_i.f_\mu$ is $x^{s_i.\mu}$.  Since $\sigma_i.f_\mu$ is a $\ttt$-eigenvector by Lemma \ref{s int}, we have $\sigma_i.f_\mu=f_{s_i.\mu}$.  If $\mu_i>\mu_{i+1}$ then by Lemma \ref{int squares} one has
\begin{align*}
\sigma_i.f_\mu&=\sigma_i^2.f_{s_i.\mu}=\frac{(z_i-z_{i+1}- c_0\pi_i)(z_i-z_{i+1}+c_0\pi_i)}{(z_i-z_{i+1})^2} f_{s_i.\mu} \\
&=\begin{cases} f_{s_i.\mu} \quad &\hbox{if $\mu_i \neq \mu_{i+1}$ mod $r$,} \\ 
\frac{(\delta-r c_0)(\delta+rc_0)}{\delta^2} f_{s_i.\mu} \quad &\hbox{if $\mu_i=\mu_{i+1}$ mod $r$.}
\end{cases}
\end{align*}
\end{proof}

\begin{corollary} \label{simple spectrum}
Suppose $\kappa=1$.  Then the $\ttt$-eigenspaces of $M(\one)$ are all one-dimensional provided $c_0 \notin \bigcup_{j=1}^{n} \frac{1}{j} \ZZ_{>0}$.
\end{corollary}
\begin{proof}
We assume there is a two dimensional $\ttt$-eigenspace and prove $c_0 \in \bigcup_{j=1}^{n} \frac{1}{j} \ZZ_{>0}$.  Let $\mu, \nu \in \ZZ_{\geq 0}^n$ be distinct and assume $\text{wt}(\mu)=\text{wt}(\nu)$; since $\text{wt}(\mu)=w.\text{wt}(\mu)$ for all $w \in S_n$ we may assume that $\mu=\mu_+$ is a partition.   Write $v=v_\nu$.  Thus $\mu_i-\mu_{i+1}=\nu_{i}-\nu_{i+1}$ mod $r$ for $1 \leq i \leq n-1$ and 
\begin{equation} \label{wts equal}
\mu_i-\nu_i=r(n-i+1-v(i))c_0
\end{equation}  for $1 \leq i \leq n$.  Let $i$ be minimal with $v(i) \neq n-i+1$.  Then $v(i)<n-i+1$ and there is some $k>i$ with $v(k)=n-i+1$.  Therefore if $c_0<0$ then
\begin{equation}
\mu_i-\nu_i=r(n-i+1-v(i))c_0<0 \quad \text{and} \quad \mu_k-\nu_k=r(n-k+1-(n-i+1)) c_0 >0,
\end{equation}  whence $\mu_i<\nu_i$ and $\mu_k>\nu_k$.  But $\mu_k>\nu_k \geq \nu_i > \mu_i$ contradicts $\mu=\mu_+$, and it follows that $c_0>0$.  Now for $1 \leq i \leq n-1$ we have
\begin{equation}
\frac{\mu_i-\nu_i-(\mu_{i+1}-\nu_{i+1})}{r}=(n-i+1-v(i)-(n-i-v(i+1))) c_0=(1+v(i+1)-v(i)) c_0
\end{equation} and the corollary follows unless $v(i+1)-v(i)=-1$ for $1 \leq i \leq n-1$.  But in that case $v=w_0$ and $\mu=\nu$, contradiction.
\end{proof}

In fact, the preceding corollary can be sharpened somewhat: provided either $p=1$ or $n$ does not divide $r$ the $\ttt$-eigenspaces are one dimensional as long as $c_0 \notin \bigcup_{j=1}^{n-1} \frac{1}{j} \ZZ_{>0}$.  We will not need this fact in this paper.

\section{Koszul resolutions of some finite dimensional $\HH$-modules} \label{L section}

We assume for the rest of the paper that $r>1$.  For $1 \leq k \leq n$ and $j \in \ZZ_{> 0}$ with $j \neq 0$ mod $r$ define (affine) hyperplanes
\begin{equation}
H_{j,k}=\left\{(c_0,d_1,\dots,d_{r/p-1}) \ | \ d_0-d_{-j}+rc_0(n-k)=j \right\} \subseteq \CC^r
\end{equation} and for $x \in \frac{1}{n} \ZZ_{>0}$,
\begin{equation}
H_x=\left\{(c_0,d_1,\dots,d_{r/p-1} \ | \ c_0=x \right\}.
\end{equation} 

The hyperplane $H_{j,1}$ was introduced (modulo different conventions for the parameters) in \cite{ChEt}, where it was called $E_j$.  Chmutova and Etingof have proved (Theorems 4.2 and 4.3 from \cite{ChEt}) that there is a finite dimensional quotient of $M(\one)$ when the parameter lies on $H_{j,1}$ for $j \neq 0$ mod $r$, and that if $p=1$ this quotient is irreducible for generic choices of the parameters.  Also, Dunkl and Opdam have proved (section 3.4 from \cite{DuOp}) that $M(\one)$ is reducible exactly if the parameter is on some $H_{j,k}$ for some positive $j \neq 0$ and $1 \leq k \leq n$ mod $r$ or $H_x$ for some $x \in \frac{1}{j} \ZZ_{>0} -\ZZ$ with $2 \leq j \leq n$.  The following theorem describes the structure of the module $M(\one)$ in the case in which $\text{gr} L(\one)$ is the quotient of the diagonal coinvariant ring predicted by Haiman.  It has the advantage of working for arbitrary divisors $p$ or $r$: this is what makes our strengthening of Vale's result possible. 

\begin{theorem} \label{submodule structure}
Suppose that $\kappa=1$, that $k \in \ZZ_{>0}$ with $k \neq 0$ mod $r$, that $(c_0,d_1,\dots,d_{r/p-1}) \in H_{k,1}$, and that the parameters do not lie on any other hyperplane $H_{l,j}$ or $H_x$ for $1 \leq l \leq n$, $j \in \ZZ_{>0}$, and $x \in \frac{1}{n} \ZZ_{>0}$.  Then the unique proper submodule of $M(\one)$ is
\begin{equation*}
\CC \left\{ f_\lambda \ | \ \hbox{$\lambda$ has at least one part of size at least $k$} \right\}.
\end{equation*} 
\end{theorem}
\begin{proof}
By Corollary \ref{simple spectrum} the $\ttt$-eigenspaces of $M(\one)$ are all one dimensional, and hence the Jack polynomials $f_\mu$ are all well-defined.  Suppose $M$ is a proper non-zero submodule of $M(\one)$.  Then $M$ contains $f_\mu$ for some $\mu \in \ZZ_{\geq 0}^n$.  By our assumption on the parameters and Lemma \ref{action lemma}, $\sigma_i.f_\mu$ is a non-zero multiple of $f_{s_i.\mu}$ whenever $\mu_i \neq \mu_{i+1}$, and it follows that $M$ also contains $f_{s_i.\mu}$ for all $1 \leq i \leq n-1$.  Hence $M$ contains $f_{\mu_-}$, where $\mu_-$ is the non-decreasing rearrangement of $\mu$.  By Lemma \ref{action lemma}, we have $\Psi.f_\mu = 0$ exactly if $\mu_n = 0$ or 
\begin{equation}
\mu_n = d_0-d_{-\mu_n}+r(v_\mu(n)-1)c_0.
\end{equation}  This last equation holds exactly if $\mu_n=k$ and $v_\mu(n)=n$, or  equivalently, exactly if $\mu_n=k$ is strictly larger than all other parts of $\mu$.  It follows that if all the parts of $\mu$ are of size less than $k$ then $M$ contains $f_0=1$, contradicting the fact that $M$ is a proper submodule.  On the other hand, if $\mu$ has at least one part of size $k$, it follows from the preceding discussion that by applying an appropriate sequence of intertwiners to $f_\mu$ we may obtain a non-zero multiple of $f_{\nu}$, where $\nu=(k,0,\ldots,0)$.  Since $f_\nu$ generates $\CC \left\{ f_\lambda \ | \ \hbox{$\lambda$ has at least one part of size at least $k$} \right\}$ as an $\HH$-module the result follows.
\end{proof}  

In Theorem 7.5 of \cite{Gri2} we generalize Theorem \ref{submodule structure} to the case of a Verma module $M(V)$ with one-dimensional $\ttt$-eigenspaces, giving a combinatorial description of the submodule structure (which can be much more intricate than the situation we study in this paper). 

\section{Diagonal coinvariants for $G(r,p,n)$} \label{diag for Grpn}

We continue to assume $r>1$.  The \emph{Coxeter number} of $G(r,p,n)$ is
\begin{equation} \label{Cox def}
h=\begin{cases}
r(n-1)+r/p \quad \hbox{if $p<r$,} \\ r(n-1) \quad \hbox{if $p=r$.}
\end{cases}
\end{equation}  This agrees with the usual definition of Coxeter number (the largest degree of a basic invariant) when $r>1$ and $p=1$ or $p=r$.  The following theorem constructs an analog for the groups $G(r,p,n)$ of the quotient of the diagonal coinvariant ring discovered by Gordon in \cite{Gor2}.  For $p<r$  a very similar theorem is proved in \cite{Val}.  Our techniques (which are conceptually very similar to those of \cite{Che}, but working directly in $\HH$) allow us to handle the case $p=r$ in the same way as $p<r$.
 
\begin{lemma} \label{free reps}
Let $m$ be a positive integer not divisible by $r$ and let $V$ be the representation $\CC \{x_1^{m},\cdots,x_n^{m} \}$.  Then $V$ is free.  If $m=h+1$ and $e_1,\dots,e_n$ are the exponents of $V$ then the multisets $\{h+1-e_i \}_{i=1}^n$ and $\{d_i \}_{i=1}^n$ are equal.
\end{lemma}
\begin{proof}
Let $\overline{m}$ and $m'$ be the integers determined by
\begin{equation}
0 \leq \overline{m}<r, \quad 0 \leq m' <r/p, \quad \overline{m}=m \ \text{mod} \ r, \quad \text{and} \quad m'=m \ \text{mod} \ r/p.
\end{equation}  Observe that the representation $\Lambda^n V$ is $\epsilon \delta^{-m}$, where $\epsilon$ and $\delta$ are the one-dimensional $G(r,1,n)$-representations determined by $\epsilon(\zeta_i^l s_ij \zeta_i^{-l})=-1$, $\epsilon(\zeta_i^l)=1$, $\delta(\zeta_i^l s_ij \zeta_i^{-l})=1$, and $\delta(\zeta_i^l)=\zeta^l$.  This is carried by the non-zero element $(x_1 \cdots x_n)^{m'} \prod_{1 \leq i<j \leq n} (x_i^r-x_j^r)$ of the ordinary coinvariant ring.    

For $1 \leq i,j \leq n$ let
\begin{equation}
f_{i,j}=x_j^{(i-1)r+\overline{m}} \quad \text{and put} \quad v_j=x_1^{r-\overline{m}+m'} \cdots x_j^{m'} \cdots x_n^{r-\overline{m}+m'}.
\end{equation}  When $p=1$ the functions $f_{i,j}$ for $1\leq j \leq n$ span a copy of $V$, and when $p>1$ $v_1,\dots,v_n$ span a copy of $V$.  One computes
\begin{equation}
\text{det}(f_{i,j})_{i,j=1}^n=(x_1 \cdots x_n)^{\overline{m}} \prod_{1 \leq i < j \leq n} (x_i^r-x_j^r)
\end{equation} and if $A$ is the matrix whose $n$th row is $v_1,v_2,\dots,v_n$ and whose $i$th row for $1 \leq i < n$ is $f_{i,1},f_{i,2},\dots,f_{i,n}$ then
\begin{equation}
\text{det}(A)=(-1)^n (x_1 \cdots x_n)^{m'} \prod_{1 \leq i < j \leq n} (x_i^r-x_j^r).
\end{equation}  It follows by Theorem 3.1 of \cite{OrSo} that $V$ is free and the exponents of $V$ are
\begin{equation}
e_i(V)=\overline{m}+(i-1)r \quad \hbox{for $1 \leq i \leq n$ if $p=1$,}
\end{equation} and
\begin{equation}
e_i(V)=\overline{m}+(i-1)r \quad \hbox{for $1 \leq i \leq n-1$ and} \quad e_n(V)=(n-1)(r-\overline{m})+nm' \quad \hbox{if $p>1$.}
\end{equation}  The degrees of $G(r,p,n)$ are $r,2r,\dots,(n-1)r, n (r/p)$ if $p<r$ and $r,2r,\dots,(n-1)r, n$ if $p=r$.  When $m=h+1$ it is straightforward to verify the last claim.
\end{proof}

\begin{theorem}
Suppose $G(r,p,n)$ acts irreducibly on $\CC^n$.  With $h$ as in \eqref{Cox def} and $V=\CC \{x_1^{h+1},\dots,x_n^{h+1} \}$, there is a $W$-equivariant quotient $L$ of the diagonal coinvariant ring $R$ of $G(r,p,n)$ such that for each $w \in W$,
\begin{equation*}
\sum \text{tr}(w,(L \otimes \Lambda^n V)_i) t^i=\frac{\text{det}(1-t^{h+1} w_V)}{\text{det}(1-t w_{\hh^*})},
\end{equation*} where $w_V$ and $w_{\hh^*}$ denote $w$ regarded as an endomorphism of $V$ and $\hh^*$, respectively.  The image of $S(\hh^*)$ in $L$ is isomorphic to the ordinary coinvariant ring.
\end{theorem}
\begin{proof}
The theorem will follow from Theorem \ref{diag coinv}, with $L=\text{gr} L(\one) \otimes \Lambda^n V^*$, once we verify its hypotheses.  Let $\mu_i \in \ZZ_{\geq 0}^n$ have an $h+1$ in the $i$th position and $0$'s elsewhere.  Then one checks that with $c_s=(h+1)/h$ for all reflections $s \in G(r,p,n)$ the hypotheses as in Theorem \ref{submodule structure} are satisfied for $k=h+1$.  Thus the radical of $M(\one)$ is generated by $\CC \{f_{\mu_1},\dots,f_{\mu_n} \}$.  Lemma \ref{action lemma} shows that as $W$-modules
\begin{equation}
\CC \{f_{\mu_1},\dots,f_{\mu_n} \} \cong V=\CC \{x_1^{h+1},\cdots,x_n^{h+1} \}.
\end{equation}  Thus the hypotheses of Lemma \ref{Koszul lemma} hold, with $k=h+1$.  The preceding lemma shows that the remaining hypotheses of Theorem \ref{diag coinv} hold.
\end{proof}

For an arbitrary irreducible complex reflection group $W$ we define the ``Coxeter'' number of $W$ to be
\begin{equation}
h=\frac{N+N^*}{n}
\end{equation} where $n$ is the dimension of the reflection reprentation of $W$, $N$ is the number of reflections in $W$ and $N^*$ is the number of reflecting hyperplanes for $W$.  This definition agrees with \eqref{Cox def} for the groups $G(r,p,n)$ whenever they are irreducible.  By a straighforward modification of \cite{BEG1} Proposition 2.3, when $c=1/h$ there is a one dimensional $\HH$-module with a BGG resolution.

\vspace{.1 in}

{\bf Question}  Is it possible that for every complex reflection group $W$ and every integer $m$ coprime to the ``Coxeter'' number $h$, when $c=(m/h)$ the $\HH$-module $L(\one)$ is $m^n$ dimensional with BGG resolution
\begin{equation}
0 \rightarrow M(\Lambda^n V) \rightarrow \cdots \rightarrow M(\Lambda^1 V) \rightarrow M(\one) \rightarrow L(\one) \rightarrow 0,
\end{equation} where $V$ is an irreducible $\CC W$-module of dimension $n$?  This question is related to Conjecture 4.3 of Bessis and Reiner \cite{BeRe}: they conjecture that for an irreducible complex reflection group of dimension $n$ that can be generated by $n$ reflections, there is a homogeneous system of parameters in each degree $\pm 1$ mod $h$ that carries either the reflection representation or its dual.  The existence of such an h.s.o.p. implies an interpretation of the $q$-Fuss/Catalan numbers as Hilbert series.  See \cite{Arm} for a survey of ``Catalan phenomena'' and non-crossing partitions.  We expect that at parameters $c_s=1+1/h$, if $W$ can be generated by $n$ reflections then $L(\one)$ gives rise to a nice quotient of the diagonal coinvariant ring; the results of Section 5 of the paper \cite{Rou} are sure to be relevant here.
\bibliographystyle{amsplain}
\def\cprime{$'$} \def\cprime{$'$}

\end{document}